\documentclass[11pt]{amsart}
\usepackage{amssymb,amsmath,amsthm}
\newtheorem{theorem}{Theorem}
\newtheorem{lemma}[theorem]{Lemma}

\newtheorem{proposition}[theorem]{Proposition}
\newtheorem{conjecture}[theorem]{Conjecture}
\newtheorem{remark}[theorem]{Remark}

\theoremstyle{remark}
\theoremstyle{definition}
\newtheorem{definition}[theorem]{Definition}
\newtheorem*{example}{Example}

\numberwithin{theorem}{section} \numberwithin{equation}{section}

\newcommand{\be}{\begin{equation}}
\newcommand{\ee}{\end{equation}}
\newcommand{\bea}{\begin{eqnarray}}
\newcommand{\eea}{\end{eqnarray}}

\newcommand{\tp}{\tilde{p}}
\newcommand{\tpp}{\tilde{p}'}

\newcommand{\ds}{\displaystyle}

\begin{document}
\title[Modular forms and graded dimensions]
{Modular forms and almost linear dependence of graded dimensions}
\author{Antun Milas}
\begin{abstract}
For every positive integral level $k$ we study arithmetic
properties of certain holomorphic modular forms associated to
modular invariant spaces spanned by graded dimensions of
$L_{\widehat{sl_2}}(k \Lambda_0)$-modules. We found a necessary
and sufficient condition for their vanishing and showed that these
modular forms resemble classical Eisenstein series
$E_{2k+2}(\tau)$. Furthermore, we derived similar results for
$\mathcal{M}(p,p')$ Virasoro minimal models, thus generalizing
some results of Mortenson, Ono and the author.
\end{abstract}

\address{Department of Mathematics and Statistics,
University at Albany (SUNY), Albany, NY 12222 }
\email{amilas@math.albany.edu}

\subjclass[2000]{11F30, 11G05, 17B67, 17B69} \maketitle

\section{Introduction}
Automorphic and modular forms are omnipresent in rational
conformal field theories (RCFTs). The most distinguished
automorphic functions that arise from a RCFT are graded dimensions
(or simply, {\em characters}) of irreducible modules of the
underlying vertex operator algebra $V$. Let us recall here that
the graded dimension of a $V$-module $M$ is defined as the graded
trace
\begin{equation} \label{char00}
{\rm ch}_M(q)={\rm tr}|_M q^{L(0)-c/24},
\end{equation}
where $L(0)$ is a Virasoro generator and $c$ is the central charge
($q=e^{2 \pi i \tau}$, for $\tau$ in the upper half-plane).

In \cite{M1},  \cite{M3} and \cite{MMO} we introduced and studied
certain remarkable modular forms associated to modular invariant
vector spaces coming from rational vertex operator algebras. In
fact, under some mild conditions, most of our results from
\cite{M1} hold for general modular invariant spaces. We briefly
recall this construction here. Let $V$ be a rational vertex
operator algebra such that modular invariance holds for
$V$-modules (see \cite{Z} for some sufficient conditions).
Consider the modular invariant vector space $\mathcal{V}$ spanned
by characters of $V$-modules. Fix an ordered basis of
$\mathcal{V}$ and let
$${F}_V(\tau)=\frac{{W'_V(\tau)}}{W_V(\tau)},$$
where $W'_V$ and $W_V$ are certain Wronskians associated to the
same ordered basis (see Section 2 for details). If $F_V$ is
nonzero we will write $\mathcal{F}_V(\tau)$ for the normalization
of $F_V$ with leading coefficient $1$. In \cite{MMO} we proved
that $\mathcal{F}_{L(c_{2,2k+1})}(\tau)$, $k \neq 6i^2-6i+1$,
resembles the classical Eisenstein series
$$E_{2k}(\tau)=1-\frac{4k}{B_{2k}} \sum_{n =1}^\infty \frac{n^{2k-1}q^n}{1-q^n}, \ k \geq 2.$$
More precisely we showed (resp. conjectured) that
$\mathcal{F}_{L(c_{2,2k+1},0)}(\tau)$ satisfy (i) and (ii) (resp.
(iii)) from the following list of properties which are known to
hold for Eisenstein series: \vskip 10mm
\begin{itemize}
\item[(i)] (Holomorphicity) For every $k \geq 2$, $E_{2k}(\tau)$
is a holomorphic modular form of weight $2k$. \item[(ii)] (Clausen
-von Staudt congruences) For prime $p=2k+1 \geq 5$, $E_{2k}(q)$ is
$p$-integral. More importantly $$ E_{2k}(q) \equiv 1 \ ({\rm mod}
\ p).$$ \item[(iii)] (Rankin-Swinnerton-Dyer \cite{R}) All zeros
of $E_{2k}(\tau)$ in the fundamental region lie on boundary arc
from $e^{2 \pi i /3}$ to $e^{\pi i/2}$ (i.e., on the $j$-line the
zeros are inside the interval $[0,1728]$).
\end{itemize}
Our purpose in this short note is to study arithmetic properties
of $F_V$ for two important families of rational vertex operator
algebras: those associated to affine Lie algebra $\widehat{sl_2}$,
denoted by $L_{\widehat{sl_2}}(k \Lambda_0)$, $k \in \mathbb{N}$
and vertex operator algebras associated to $\mathcal{M}(p,p')$
Virasoro minimal models, denoted by $L(c_{p,p'},0)$, where
$$c_{p,p'}=1-\frac{6(p-p')^2}{pp'}, \ (p,p')=1, \ p \geq 2,
p' \geq 2.$$ For more about these vertex operator algebras see
\cite{LL}. To simplify the notation we shall write $L(k
\Lambda_0)$ instead of $L_{\widehat{sl_2}}(k \Lambda_0)$.

%In our joint work with
%Mortenson and Ono we studied the family the case
%$V=L(c_{2,2k+1},0)$, where $k \geq 2$ \cite{MMO}.
As we already noticed in \cite{MMO} and \cite{M3} a nontrivial
task is to determine precisely the vanishing condition for
${F}_{V}(\tau)$. Our first result establishes this for $F_{L(k
\Lambda_0)}$.
\begin{theorem} \label{one}
For every $k \geq 1$, $F_{L(k \Lambda_0)}$ is a holomorphic
modular form. Moreover, $F_{L(k \Lambda_0)}(\tau)=0$ if and only
if $k=2i^2-2$ for some $i \geq 1$.
\end{theorem}
\noindent Similarly,
\begin{theorem} \label{onep}
For every $p$ and $p'$, $F_{L(c_{p,p'},0)}$ is a holomorphic
modular form. We have $F_{L(c_{p,p'},0)}(\tau)= 0$ if and only if
$(p,p')=(2\tilde{p}^2,3 \tilde{p}'^2)$.
\end{theorem}
The previous theorem has been established for $L(c_{2,2k+1},0)$,
$k \geq 2$ in \cite{MMO}. In the process of proving Theorems
\ref{one} and \ref{onep} an important ingredient is an application
of $q$-series identities among irreducible characters of the form
\be \label{near-miss} \sum_{i=1}^{s} m_i {\rm ch}_{M_i}(\tau)=C
\neq 0, \ee where $M_i$ are irreducible $V$-modules, $C$ is a
constant and $m_i \in \mathbb{Z}$. The equation (\ref{near-miss})
is what we refer to as an "almost linear dependence" relation
among graded dimensions ${\rm ch}_{M_i}$.

In order to prove (ii) for $\mathcal{F}_V$, we have to prove the
$p$-integrality first. This is a nontrivial problem because in
some cases the leading term of $F_V$ is divisible by $p$.
\begin{theorem} \label{two}
For every prime $p=2k+3 \geq 5$, $\mathcal{F}_{L(k
\Lambda_0)}(\tau)$ is $p$-integral.
\end{theorem}
This result is far from being obvious, simply because the leading
coefficient of $F_{L(k \Lambda_0)}$ is divisible by $p=2k+3$ for
{\em every} $p$. Based on overwhelming numerical evidence, Theorem
\ref{one} and Theorem \ref{two} we conjecture:
\begin{conjecture} \label{myfirst}
For $k \neq 2i^2-2$, $i \geq 1$ and prime $p=2k+3 \geq 5$ we have
\be \label{hasse} \mathcal{F}_{L(k \Lambda_0)}(q) \equiv 1 \ ({\rm
mod} \ p). \ee
 Moreover, the zeros of $\mathcal{F}_{L(k \Lambda_0)}$,
in the fundamental domain, lie on boundary arc from $e^{2 \pi i
/3}$ to $e^{\pi i/2}$.
\end{conjecture}
For $\mathcal{M}(p,p')$ minimal models the property (iii) fails to
hold in general (e.g., $\mathcal{F}_{L(c_{3,7},0)}$). Because of
complexity of computation we do not have enough numerical evidence
which would support a precise conjecture for
$\mathcal{F}_{L(c_{p,p'},0)}$. Nevertheless, we did observe that
properties (ii) and (iii) seem to hold for unitary minimal models
$\mathcal{M}(m+2,m+3)$.

\noindent {\bf Acknowledgments.} This paper is dedicated to Jim
Lepowsky and Robert L. Wilson on the occasion of their 60th
birthday.
%They made an
%important contribution in discovering a link between classical
%$q$-series identities and infinite-dimensional Lie algebras.

We would like to thank K. Ono for discussion on related topics, to
E. Mortenson for providing us with enough computational evidence
that backed Conjecture \ref{myfirst}. Finally, we thank E. Mukhin
for a correspondence.

\section{Irreducible characters of $L(k \Lambda_0)$--modules and $L(c_{p,p
'},0)$-modules.}

In this section we are following \cite{K}.
Let $k$ be a positive integer (level). For every $i=1,...,k+1$,
let $L((k-i+1)\Lambda_0+(i-1)\Lambda_1)$ be the integrable highest
weight $\widehat{sl_2}$-module of level $k$ with highest weight
$(k-i+1)\Lambda_0+(i-1)\Lambda_1$. For simplicity we shall write
$L(k,i-1)$ instead of $L((k-i+1)\Lambda_0+(i-1)\Lambda_1)$. It is
known that $L(k \Lambda_0)$ has a natural vertex operator algebra
structure, but we will not use this fact in the rest of the text.
The conformal vector $\omega_k \in L(k \Lambda_0)$ is given by
$$\omega_k=\frac{1}{2(k+2)}\left(\frac{h(-1)^2{\bf
1}}{2}+e(-1)f(-1){\bf 1}+f(-1)e(-1){\bf 1}\right),$$ where ${\bf
1}$ is the vacuum vector and $\{e,f,h \}$ is the standard basis of
the finite-dimensional simple Lie algebra $\frak{sl}_2$. The
central charge of $L(k \Lambda_0)$ is given by
\begin{equation}
c_k = \frac{3k}{k+2}, \nonumber
\end{equation}
and the lowest conformal weight of $L(k,i-1)$ is
\begin{equation} h_{k,i}
= \frac{i^2-1}{4(k+2)}, \nonumber
\end{equation}
where $i=1,...,k+1$. Thus,
\begin{equation}
h_{k,i}-\frac{c_k}{24}=\frac{i^2}{4(k+2)}-\frac{1}{8}.
\end{equation}

The homogeneous specialization
%(i.e., $e^{-\alpha_0}=q$ and
%$e^{-\alpha_1}=1$)
in the Weyl-Kac character formula yields (cf. \cite{K}): For every
positive integer $k$, and $i=1,...,k+1$, we have
\begin{equation} \label{irr-char}
{\rm ch}_{k,i}(q):={\rm
tr}|_{L(k,i-1)} q^{L(0)-c/24}=\frac{\ds{\sum_{\tiny{\begin{array}{cc} n \in \mathbb{Z} \\
n \equiv \ i \ {\rm mod} \ 2(k+2)
\end{array}}} n q^{n^2/4(k+2)}}}{\eta(\tau)^3},
\end{equation}
where $\eta(\tau)=q^{1/24} \ds{\prod_{i=1}^\infty(1-q^i)}$ is the
Dedekind eta-functions.
%Again, for simplicity, we shall write
%${\rm ch}_{k,i}(q)$ instead of ${\rm ch}_{L(k,i-1)}(q)$.

%There is a distinguished family or rational vertex operator
%algebras associated to certain lowest weight representation of the
%Virasoro algebra (see \cite{LL} for instance).
Let $p \geq 2$ and $p' \geq 2$ be two relatively prime positive
integers. We have the following well-known formula for the graded
dimensions of $L(c_{p,p'},0)$-modules (cf. \cite{RC}). \be
\label{vir-char} {\rm ch}_{p,p'}^{r,s}(q):={\rm
tr}|_{L(c_{p,p'},h_{p,p'}^{r,s})} q^{L(0)-c_{p,p'}/24} =\frac{
\ds{\sum_{n \in \mathbb{Z}} q^{\frac{(2npp'+p'r-ps)^2}{4p
p'}}-q^{\frac{(2npp'+p'r+ps)^2}{4pp'}}}}{\eta(\tau)}, \ee where $1
\leq r \leq p-1$ and $1 \leq s \leq p'-1$.

\section{Modular forms and $SL(2,\mathbb{Z})$--modules}

Form now on we will be using the Ramanujan's derivative $(q
\frac{d}{dq})$. Let $\{ f_1,....,f_m \}$ be a basis of the modular
invariant space spanned by characters of a rational vertex
operator algebra $V$. Let also $W_V(\tau)=W_{(q
\frac{d}{dq})}(f_1,...,f_m)$, where $W_{(q \frac{d}{dq})}( \cdot
)$ denotes the Wronskian determinant with respect to $q
\frac{d}{dq}$. We will denote by $\mathcal{W}_V$ a multiple of
$W_{(q \frac{d}{dq})}(f_1,...,f_m)$ with the leading coefficient
$1$. This normalization is important in the theory of mod $p$
modular forms \cite{Se}. Similarly, we will denote by $W'_V=W_{(q
\frac{d}{dq})}(f'_1,...,f'_m)$ the Wronskian of derivatives of
$f_i$, and by $\mathcal{W}'_V$ its normalization (if nonzero) with
the leading coefficient being $1$. Now, $W_V$ and $W'_V$ are
determined only up to a nonzero constant, while $\mathcal{W}_V$,
$\mathcal{W}'_V$,
$$F_V(\tau):=\frac{W'_V(\tau)}{W_V(\tau)},$$
and $\mathcal{F}_V:=\frac{\mathcal{W}'_V}{\mathcal{W}}$ do not
depend on the choice of a basis of $\mathcal{V}$.
%If the span of $f_1,...,f_k$ is a $SL(2, \mathbb{Z})$--module
%under the modular transformation it is not hard to see that $W_{(q
%\frac{d}{dq})}(f_1,...,f_k)$ is an automorphic form (i.e., a
%modular form with a character).
If $W_{(q \frac{d}{dq})}(f_1,...,f_m)$ has no zeros in the upper
half-plane, the quotient $F_V(\tau)$ is a holomorphic modular form
of weight $2m$ (see \cite{M1}, \cite{M3}).
% From now on we will Let
%$V$ is a finite-dimensional modular invariant space spanned by
%irreducible characters of $V$-modules.
%Unlike $W_V$ and $W'_V$,
%$$F_V(\tau):=\frac{W'_V(\tau)}{W_V(\tau)}$$
%does not depend on the choice of the basis of $V$.

%We can make everything more explicit in the case of irreducible
%characters of $\hat{\frak{sl}_2}$.
The following result is from \cite{M1}
\begin{theorem} \label{ded-1}
For every $k \geq 1$,
$$\mathcal{W}_{L(k \Lambda_0)}(q)=\eta(q)^{2k(k+1)}.$$
Consequently, ${F}_{L(k\Lambda_0)}$ is a holomorphic modular form
of weight $2k+2$ (possibly zero).
\end{theorem}
Similarly, we have a result from \cite{M0} (see also \cite{M1} for
a different proof).
\begin{theorem}
For every $p$ and $p'$ as above
$$\mathcal{W}_{L(c_{p,p'},0)}(q)=\eta(q)^{\frac{(p-1)(p'-1)(pp'-p-p'-1)}{2}}$$
Thus, ${F}_{L(c_{p,p'},0)}$ is a holomorphic modular form of
weight $\frac{(p-1)(p'-1)(pp'-p-p'-1)}{4}$ (possibly zero).
\end{theorem}

\section{Vanishing results for ${F}_{L(k \Lambda_0)}$ and triangular numbers}

Let us recall a classical Jacobi's $q$-series identity \be
\label{jacobi3} \prod_{n=1}^\infty (1-q^n)^3=\sum_{n=0}^\infty
(-1)^n (2n+1)q^{n(n+1)/2}. \ee This $q$-series identity has a nice
representation theoretic interpretation.
%It is easy to see that $h_{k,i}=\frac{c_k}{24}$ if and only if
%$i^2=2(k+2)$ or equivalently $k=2i^2-2$, $i \geq 2$. Moreover,
%$$h_{k,i}-\frac{c_k}{24} \in \mathbb{N},$$
%if and only if $k=2i^2-2$ and $j(2i+1)$, $j=1,...,i$.
We will need the following lemma.
\begin{lemma} \label{trivial}
For every positive integer $i \geq 2$,
\begin{eqnarray}
&& \prod_{n=1}^\infty (1-q^n)^3=\sum_{j=0}^{i-1} \sum_{m=0}^\infty
(4mi+2j+1)(-1)^j
q^{(2mi+j)(2mi+j+1)/2} \nonumber \\
&&+ \sum_{j=0}^{i-1} \sum_{m=1}^\infty (4mi-2j-1)(-1)^{j+1}
q^{(2mi-j)(2mi-j-1)/2}.
\end{eqnarray}
\end{lemma}
\begin{proof}
For every $n \in \mathbb{N}$ let $n \equiv j \ ({\rm mod} \ 2i)$,
where $j  \in \{-i,...,-1,0,1,...,i-1 \}$. Now, we apply
(\ref{jacobi3}) and rewrite the sum $\ds{\sum_{n=0}^\infty (-1)^n
(2n+1)q^{n(n+1)/2}}$ modulo $2i$ and sum over $j$ as above.
\end{proof}
\noindent Firstly, we classify all $k$, such that for some $i$,
${\rm ch}_{k,i}(q)$ has a nonzero leading constant term and
consequently only positive integer powers. It is easy to see that
this is equivalent to
$$h_{k,i}=\frac{c_k}{24},$$ which holds if and only if
$2i^2=(k+2)$, or equivalently $k=2i^2-2$, for some $i \geq 2$.
Moreover, for $k=2i^2-2$
$$h_{2i^2-2,m}-\frac{c_{2i^2-2}}{24} \in \mathbb{N},$$
if and only if $m=i(2j+1)$, $j=0,...,i$.

\begin{theorem} We have $\mathcal{F}_{L(k \Lambda_0)}(q)=0$ if and only if
$k=2i^2-2$ for some $i \geq 1$.
\end{theorem}
\begin{proof}
We rewrite the irreducible characters (\ref{irr-char}) in the
following form:
\begin{equation} \label{irr-char1}
{\rm ch}_{k,i}(q)=\frac{\ds{\sum_{m \in \mathbb{Z}}}
(i+2m(k+2))q^{(i+2m(k+2))^2/4(k+2)}}{\eta(\tau)^3}.
\end{equation}
Now, let us use the fact that
$${\rm ch}_{k,i}(q) \in a_0+ q \mathbb{Z}_{\geq 0}[[q]], \ \ a_0 \neq 0$$
if and only if $k=2i^2-2$ and for such $k$, for $j=0,..,i-1$, we
have
$${\rm ch}_{(2i^2-2),(2j+1)i}(q) \in \mathbb{Z}_{\geq 0}[[q]].$$
Now, for $j=0,...,i-1$
$${\rm ch}_{(2i^2-2),(i(2j+1))}(q)=\frac{\ds{\sum_{m \in \mathbb{Z}}
i(4mi+2j+1)q^{\frac{(j+2mi)(j+2mi+1)}{2}}}}{(q;q)_\infty^3}.$$ Let
us rewrite the numerator. Clearly,
\begin{eqnarray} \label{split}
&& \sum_{m \in \mathbb{Z}}
i(4mi+2j+1)q^{\frac{(j+2mi)(j+2mi+1)}{2}}=i \sum_{m=
0}^\infty(4mi+2j+1)q^{\frac{(j+2mi)(j+2mi+1)}{2}}  \nonumber \\ &&
+ i \sum_{m=0}^\infty
(-1)(4mi+4i-2j-1)q^{\frac{(2mi+2i-j)(2mi+2i-j-1)}{2}}.
\end{eqnarray}
Now,
\begin{eqnarray}
&& \sum_{j=0}^{i-1} (-1)^j \sum_{m \in \mathbb{Z}}
i(4mi+2j+1)q^{\frac{(j+2mi)(j+2mi+1)}{2}}=\sum_{j=0}^{i-1}
\sum_{m=0}^\infty (-1)^j (4mi+2j+1)q^{\frac{(j+2mi)(j+2mi+1)}{2}}
\nonumber \\
&& + \sum_{j=0}^{i-1} \sum_{m=0}^\infty (-1)^{j+1}
(4mi+4i-2j-1)q^{\frac{(2mi+2i-j)(2mi+2i-j-1)}{2}} \nonumber \\
&& = \sum_{j=0}^{i-1} \sum_{m=0}^\infty (-1)^j
(4mi+2j+1)q^{\frac{(j+2mi)(j+2mi+1)}{2}}+ \sum_{m=1}^\infty
(-1)^{j+1}(4mi-2j-1)q^{\frac{(2mi-j)(2mi-j-1)}{2}}. \nonumber
\end{eqnarray}
By Lemma \ref{trivial} and formula (\ref{split}) we have
$$ \sum_{j=0}^{i-1} (-1)^j {\rm ch}_{2i^2-2,j(2i+1)}(q)=i.$$
Therefore ${W}'_{L((2i^2-2)\Lambda_0)}(q)$ is zero.

Now, we prove the converse. If ${W}'_{L(k \lambda_0)}$ is zero,
then there exist ${\rm ch}_{k,i}(q)$, $i \in I \subseteq
\{0,1,...,k\}$ such that $$\sum_{i \in I} m_i {\rm
ch}'_{k,i}(\tau)=0.$$ On the other hand, ${\rm ch}_{k,i}(q)$, $i
\in I$ are linearly independent, so \be \label{sufficient} \sum_{i
\in I} m_i {\rm ch}_{k,i}(\tau)=C \neq 0, \ee for some constant
$C$. The previous equation implies that there exists a subset $J
\subseteq I$, such that ${\rm ch}_{k,j}(q)$, $j \in J$ admits only
positive integer powers of $q$. From (\ref{sufficient}) it follows
that
$${\rm ord}_{i \infty}({\rm ch}_{k,j})=0,$$
for some $j$. Therefore, $k=2i^2-2$ for some $i$.
%But we already
%know that this is a sufficient condition for vanishing.
\end{proof}

\section{Vanishing results for $F_{L(c_{p,p'},0)}$ and pentagonal numbers}

In this section we will give necessary and sufficient conditions
for vanishing of ${F}_{L(c_{p,p'},0)}$. We will use approach from
the previous section and \cite{MMO}. The following lemma gives a
necessary condition for the vanishing.
\begin{lemma} \label{one-lemma} If $F_{L(c_{p,p'},0)}(q)=0$, then $pp'=6m^2$,
for some $m \in \mathbb{N}$.
\end{lemma}
\begin{proof} As in the previous section, it is not hard to see that the vanishing
of $F_{L(c_{p,p'},0)}$ implies that there exists a pair $(r,s)$,
such that ${\rm ord}_{i \infty}({\rm ch}_{p,p'}^{r,s})=0$. Now the
equation $h_{p,p'}^{r,s}-\frac{c_{p,p'}}{24}=0$ implies the wanted
condition.
\end{proof}
The previous lemma gives the following four
\footnote{$\mathcal{M}(p,p')$ and $\mathcal{M}(p',p)$ are the same
minimal models, so there are actually only two cases to consider.}
possibilities on $p$ and $p'$: \bea \label{1234} && p=\tilde{p}^2,
\ \ \ p'=6\tpp^2, \nonumber
\\ && p=6\tilde{p}^2, \ \ p'=\tpp^2, \nonumber \\ &&
p=3\tilde{p}^2, \ \ p'=2\tpp^2 \nonumber \\ && p=2\tilde{p}^2, \ \
p'=3\tpp^2. \nonumber \eea We will rule out the first case a
little bit later. Now, let $p=2 \tp^2$ and $p'=3 \tpp^2$. We shall
classify all ${\rm ch}_{p,p'}^{r,s}(q)$ with positive integer
powers of $q$. A simple computation shows that ${\rm
ch}_{p,p'}^{r,s}(q)$ has positive integer powers of $q$ if and
only if
$$6(p'r-sp)^2=pp'(24k+1), \ \ {\rm for} \ {\rm some} \ k \geq 0.$$
Now, in Lemma \ref{one-lemma} we already established that
$pp'=6m^2$. Thus, $(p'r-sp)^2=m^2(24k+1)$ implies that $24k+1$ is
a perfect square. Equivalently, $k=\frac{3l^2+l}{2}$ for some $l
\in \mathbb{Z}$. In this case \be \label{pentagonal} {\rm ord}_{i
\infty} ({\rm ch}_{p,p'}^{r,s}(q))=\frac{3l^2+l}{2}. \ee Now we
have to classify all $(r,s)$-pairs satisfying (\ref{pentagonal}).
The equation \be \label{cons} (p'r-sp)^2=\tp^2 \tpp^2(24k+1) \ee
adds an additional constraint on $r$ and $s$. Let us consider the
case $p=2 \tp^2$, $p'=3 \tpp^2$ for $\tp$ and $\tpp$ relatively
prime. It follows then that $\tp |r$ and $\tpp |s$, so that
$s=\tpp s'$ and $r=\tp r'$, for some $r'$ and $s'$ such that $1
\leq r' \leq 2 \tp-1$ and $1 \leq s' \leq 3 \tpp-1$. If we
substitute everything in (\ref{cons}) and factor $\tp \tpp$ we
obtain $(3r'\tpp-2 s' \tp)^2=(24k+1)$, which holds if and only if
$(3r' \tpp -2 s'\tp) \equiv \pm 1 \ ({\rm mod} \ 6)$. Clearly,
$\tp \equiv \pm 1 \ ({\rm mod} \ 3)$ and $\tpp \equiv 1 \ ({\rm
mod} \ 2)$.
%Thus $3r'\tpp-2 s'
%\tp \equiv \pm 1 \ ({\rm mod} \ 6)$ if and only if $s' \equiv \pm
%1 \ ({\rm mod} \ 3)$ and $r' \equiv \pm 1 \ ({\rm mod} \ 2)$.
Now if $s' \equiv -1 \ ({\rm mod} \ 3)$, then because of the
symmetry ${\rm ch}_{2\tp^2,3\tpp^2}^{\tp r',\tpp s'}(q)={\rm
ch}_{2\tp^2,3\tpp^2} ^{2 \tp^2-\tp r',3 \tpp^2-\tpp s'}(q)$, it
suffices to consider only $r'$ and $s'$ such that $r'$ is odd and
$s' \equiv 1 \ ({\rm mod} \ 3)$. This completes the classification
of ${\rm ch}_{2\tp^2,3\tpp^2}^{r,s}(q)$ with positive integer
powers.

%As an example consider $\mathcal{M}(6,25)$ minimal models. The
%pair $p=6$, $p'=25$ satisfies the condition in Lemma
%\ref{one-lemma}, but $\mathcal{F}_{L(c_{6,25},0)} \neq 0$.

The following result is a generalization of a series of identities
for $\mathcal{M}(2,3(2k-1)^2)$ Virasoro minimal model characters
from \cite{MMO}. Essentially the same formula was obtained
recently by Mukhin \cite{Mu} by using similar methods, but
apparently motivated by a different circle of ideas (cf.
\cite{BF}).
\begin{proposition} \label{me-mukhin-prop}
The following $q$-series identity holds among characters of
$\mathcal{M}(2\tp^2,3 \tpp^2)$ minimal models \be
\label{me-mukhin} \ds{\sum_{{r'}=1, r' \equiv 1 \ ({\rm mod} \ 2)
}^{2(\tp-1)+1} \sum_{{s'}=1, s' \equiv 1 \ ({\rm mod} \
3)}^{3(\tpp-1)+1} (-1)^{\frac{3r'\tp-2s'\tpp+1}{2}} {\rm ch}_{2
\tp^2,3\tpp^2}^{r' \tp ,s' \tpp}(q)=1.} \ee
\end{proposition}
\begin{proof} The idea is the same as in the previous chapter and \cite{MMO} so we will omit
some details. Let us recall Euler's Pentagonal Number Theorem
\cite{A}: \be \label{euler} \eta(\tau)=\sum_{k \in \mathbb{Z}}
(-1)^k q^{\frac{(6k+1)^2}{24}}. \ee
 Now, from (\ref{vir-char}) it
follows that \be \label{char-mod} {\rm ch}_{2\tp^2,3\tpp^2}^{r'
\tp,s' \tpp}(q)=\eta^{-1}(\tau) \sum_{m \in \mathbb{Z}} \left(
q^{\frac{(12 \tp \tpp m+3r' \tpp -2 s'\tp)^2}{24}}-q^{\frac{(12
\tp \tpp m+3r' \tpp + 2 s'\tp)^2}{24}} \right). \ee Notice that
all powers of $q$ in the numerator of (\ref{char-mod}) are
pentagonal numbers shifted by $\frac{1}{24}$. Now, for $r' \equiv
1 \ ({\rm mod} \ 2)$ and $s' \equiv 1 \ ({\rm mod} \ 3)$, the
numbers $3r' \tpp-2s'\tp$ and $3r' \tpp+2 s' \tp$ are all distinct
($2\tp \tpp$ values in total) and congruent to $\pm 1$ modulo $6$.
Thus {\small \bea \label{sets} && \biggl\{ \frac{(6(2 \tp \tpp
m)+3r' \tpp -2 s'\tp)^2}{24} : m \in \mathbb{Z}, 1 \leq r' \leq 2
\tp-1, r' \equiv 1 \ ({\rm mod} \ 2); 1 \leq s' \leq 3 \tpp-1, s'
\equiv 1 \ ({\rm mod} \ 3) \biggr\} \nonumber
\\ && \bigcup \biggl\{ \frac{(6(2 \tp \tpp m)+3r' \tpp + 2 s'\tp)^2}{24}
: m \in \mathbb{Z}, 1 \leq r' \leq 2 \tp-1, r' \equiv 1 \ ({\rm
mod} \ 2); 1 \leq s'
\leq 3 \tpp-1, s' \equiv 1 \ ({\rm mod} \ 3) \biggr\} \nonumber \\
&&=\biggl\{ \frac{(6(2 \tp \tpp m+i)+1)^2}{24} : \ -p p' < i \leq
pp' , m \in \mathbb{Z} \}=\biggl\{ \frac{(6k+1)^2}{24}, k \in
\mathbb{Z} \biggr\}. \nonumber \eea } Finally, the formula
(\ref{char-mod}) gives
$$\sum_{{\tiny \begin{array}{c} r'=1 \\ r' \equiv 1 \ ({\rm mod} \ 2 ) \end{array}}}^{2(\tp-1)+1} \sum_{{\tiny \begin{array}{c} s'=1 \\
s' \equiv 1 \ ({\rm mod} \ 3) \end{array} }}^{3 (\tpp-1)+1}
(-1)^{\frac{ 3 r' \tpp- 2 s' \tp+1}{2}} {\rm
ch}_{2\tp^2,3\tpp^2}^{r' \tp,s' \tpp}(q) =\frac{\ds{\sum_{ t=-\tp
\tpp+1}^{\tp \tpp} \sum_{m \in \mathbb{Z}} (-1)^{t}
q^{\frac{(6(2\tp \tpp m+t)+1)^2}{2}}}}{\eta(\tau)}.$$ The last
expression equals $1$ by the Euler's identity (\ref{euler}).
\end{proof}
%$n \equiv
%k \ ({\rm mod} \ 2pp')$, $k \in \{- \tp \tpp,-\tp \tpp,...,\tp
%\tpp-1 \}$. It is easy to see that for $r'$ and $s'$ as abo

\begin{lemma} \label{two-lemma}
If $p=\tp^2$ and $p'=6 \tpp^2$, then $F_{L(c_{p,p'},0)} \neq 0$.
\end{lemma}
\begin{proof}
It suffices to consider the case $p=\tp^2$ and $p'=6 \tpp^2$. Then
$\tp \equiv \pm 1 \ ({\rm mod} \ 6)$. Here the $q$-series ${\rm
ch}_{\tp^2,6\tpp^2}^{r,s}(q)$ with integer powers are given by \be
\label{char-mod-6} {\rm ch}_{\tp^2,6\tpp^2}^{r' \tp,s'
\tpp}(q)=\eta^{-1}(\tau) \sum_{m \in \mathbb{Z}} \left(
q^{\frac{(12 \tp \tpp m+6r' \tpp -s'\tp)^2}{24}}-q^{\frac{(12 \tp
\tpp m+6r' \tpp + s'\tp)^2}{24}} \right), \ee where $1 \leq r'
\leq \tp-1$ and $1 \leq s' \leq 6\tpp-1$. Now, we have to show
that no linear combination of numerators in (\ref{char-mod-6})
equals $\eta(\tau)$. Thus, it suffices to show that there exists
$k \in \mathbb{Z}$,  such that $q^{\frac{(6k+1)^2}{24}}$ does not
appear in the numerator of (\ref{char-mod-6}) for any admissible
choices of $r'$ and $s'$. Because of $(p,p')=1$, we have  $\tp
\equiv \pm 1 \ ({\rm mod} \ 6)$. But equations $12 \tp \tpp m+
r'\tpp-6 s'\tp=\tp$ and $12 \tp \tpp m+ r'\tpp+6 s'\tp=\tp$ have
no solutions in the given range for $r'$ and $s'$, thus it is
impossible to find a linear combination of characters ${\rm
ch}_{\tp^2,6 \tpp^2}^{r,s}$ which equals to $1$.
\end{proof}

\begin{theorem}
The modular form $F_{L(c_{p,p'},0)}$ is zero if and only if $p=2
\tp^2$, $p'=3\tpp^2$.
\end{theorem}
\begin{proof}
If $p =2 \tp^2$ and $p'=3 \tpp^2$ the result follows from Lemma
\ref{one-lemma}, the discussion after Lemma \ref{one-lemma},
Proposition \ref{me-mukhin-prop} and Lemma \ref{two-lemma}.
\end{proof}

\begin{remark}
{\em For fixed $p$ and $p'$, the $q$-series identities obtained in
Proposition \ref{me-mukhin} is the only almost-linear dependence
relation of the form \be \label{alm-lin}
 \sum_{i \in I} m_i {\rm ch}_{p,p'}^{r,s}(q)=1, \ m_i \neq
0. \ee Otherwise the irreducible characters would be linearly
dependent, which is false.}
\end{remark}

\section{$p$-integrality of $\mathcal{F}_{L(k \Lambda_0)}$, $p=2k+3$}

As usual, we say that $a \in \mathbb{Q}$ is $p$-integral if
${v}_{p}(a) \geq 0$. For two formal $q$-series $$A(q)=\sum_{n \in
\mathbb{C}} a(n)q^n \in \mathbb{Q} \{q \}, \ \ B(q)=\sum_{n \in
\mathbb{C}} b(n) q^n \in \mathbb{Q} \{q \}$$ with $p$-integral
coefficients we will write
$$A(q) \equiv B(q) \ ({\rm
mod} \ p),$$ if
$$a(n) \equiv b(n) \ ({\rm mod} \ p), \ \ {\rm for \ \  every} \ n.$$
\begin{lemma}
For $k \geq 2$ and $p=2k+3$ prime, the $q$-series ${W}_{L(k
\Lambda_0)}(q)$, $\mathcal{W}_{L(k \Lambda_0)}(q)$ and $W'_{L(k
\Lambda_0)}(q)$ are $p$-integral.
\end{lemma}
\begin{proof}
The infinite product $\mathcal{W}_{L(k \Lambda_0)}(\tau)$ (cf.
Theorem \ref{ded-1}) is clearly $p$-integral. Furthermore $4(k+2)$
is not divisible by $p$, so it follows that ${W}_{L(k
\Lambda_0)}(\tau)$ and ${W}'_{L(k \Lambda_0)}$ are $p$-integral as
well.
\end{proof}
\noindent It is not clear at all that the $q$-series
$\mathcal{W}'_{L(k \Lambda_0)}(q)$ is $p$-integral.

\noindent For a formal series $F(y) \in \mathbb{Q}\{y,q\}$, let
$$[y^k]F(y):={\rm Coeff}_{y^k} F(y) \in \mathbb{Q}\{q\}.$$
\begin{definition}
Let $F(q)$ be a (formal) $q$-series. Then the $s$-th moment of
$F(q)$ is defined as
$$\frac{\left(q \frac{d}{dq} \right)^s F(q)}{F(q)}.$$
\end{definition}

\begin{proposition} \label{jacmodp}
\begin{itemize}
\item[(a)] For every $k \geq 2$, and $p=2m+1$ prime we have the
following congruence
\begin{equation} \label{lhs}
m! [y^m] \left(\frac{\eta(qe^y)}{\eta(q)}\right)^3 \equiv 2^{-3m}
\ ({\rm mod} \ p).
\end{equation}
\item[(b)] For every $m$ the left hand side in (\ref{lhs}) is a
quasimodular form.
\end{itemize}
\end{proposition}
\begin{proof}
By using Jacobi's formula we have
$$\eta^3(q e^y)=\sum_{n=0}^\infty (-1)^n (2n+1) e^{y(2n+1)^2/8}
q^{(2n+1)^2/8}.$$ Thus,
\begin{eqnarray} \label{gm}
m![y^m]\eta(qe^y)^3 &=& \sum_{n=0}^\infty (-1)^n
\frac{(2n+1)^{2m+1}}{8^m} q^{(2n+1)^2/8} \nonumber \\
&\equiv& 2^{-3m} \sum_{n=0}^\infty (-1)^n (2n+1) q^{(2n+1)^2/8} \
({\rm mod} \ p).
\end{eqnarray}
 Now,
$$m![y^m]\left(\frac{\eta(e^y q)}{\eta(q)}\right)^3 \equiv 2^{-3m} \
({\rm mod} \ p).$$ Part (b) follows from the fact that the
logarithmic derivative of $\eta(\tau)$ is (up to a non-zero
constant) the quasimodular form $E_2(\tau)$ and the ring of
quasimodular forms is closed under the differentiation $\left(q
\frac{d}{dq}\right)$ (cf. \cite{KZ2}).
\end{proof}
\begin{example}
Let $E_{2m,3}(\tau)$ denote the normalization of the series on the
left hand side of (\ref{gm}) such that the leading coefficient in
the $q$-expansion is one, then
\begin{eqnarray}
E_{2,3}(\tau)&=& E_2(\tau), \nonumber \\
E_{4,3}(\tau)&=&\frac{5}{3}E_2^2(\tau)-\frac{2}{3}E_4(\tau),
\nonumber \\
E_{6,3}(\tau)&=&
\frac{35}{9}E_2(\tau)^3-\frac{14}{3}E_2(\tau)E_4(\tau)+\frac{16}{9}E_6(\tau).
\nonumber
\end{eqnarray}
\end{example}
Let us denote the numerator in (\ref{irr-char1}) by
$\theta_{k,i}(\tau)$, so that
\begin{equation} \label{irr-char2}
\theta_{k,i}(\tau)=\sum_{m \in \mathbb{Z}}
(i+2m(k+2))q^{(i+2m(k+2))^2/4(k+2)}.
\end{equation}
Let $p=2k+3$ be prime. Clearly, $GCD(4(k+2),p)=1$.
\begin{lemma} \label{k1st}
For every $i=1,...,k+1$, we have
\begin{equation} \label{congruence}
\theta^{(k+1)}_{k,i}(\tau) \equiv (4(k+2))^{-k-1}
\theta_{k,i}(\tau) \ ({\rm mod} \ p).
\end{equation}
\end{lemma}

\noindent Now we are ready to prove
\begin{theorem} \label{ss}
For every $k \geq 1$ such that  $p=2k+3$ is prime,
$\mathcal{W}'_{L(k \Lambda_0)}$ (and henceforth $\mathcal{F}_{L(k
\Lambda_0)}$) is $p$-integral.
\end{theorem}
\begin{proof}
It suffices to prove
\begin{equation} \label{a}
\frac{W'_k(\tau)}{W_k(\tau)} \equiv 0 \ ({\rm mod} \ p)
\end{equation}
and
\begin{equation} \label{b}
v_p(a_0)=1, \ \ {\rm where} \ \ \frac{W'_k}{W_k}=a_0+a_1q+\cdots.
\end{equation}
We first prove (\ref{b}). Notice that
$$p \nmid h_{k,i}-\frac{c_k}{24},$$
for $i=1,...,k$, but $p$ divides $$
h_{k,k+1}-\frac{c_k}{24}=\frac{(2k+3)k}{8(k+2)}.$$ Thus, the
leading coefficient $a_0$ in the $q$-expansion of $W'_{L(k
\Lambda_0)}$ is divisible by $p$. In fact, $$a_0=\prod_{i=1}^{k+1}
\left(h_{k,i}-\frac{c_k}{24}\right),$$ so that $v_p(a_0)=1$. It is
easy to see that the leading coefficient in the $q$-expansion of
$W_{L(k \Lambda_0)}$ is (up to a sign) $\ds{\prod_{1 \leq m < n
\leq k+1} \frac{m^2-n^2}{4(k+2)}}$, which is not divisible by
$p=2k+3$. This proves (\ref{b}).

%Notice first that $2(2i^2-2)+3=4i^2-1=(2i-1)(2i+1)$, thus
%$\mathcal{W}_k(\tau) \neq 0$ for $\frac{k-3}{2}=p$, $p$ prime.
%$$\mathcal{F}_k=\frac{\mathcal{W}'_k(\tau)}{\mathcal{W}_k(\tau)} \equiv 1  \ ({\rm mod} \ p).$$
%Let us recall that $\mathcal{W}_k(\tau)$ and
%$\mathcal{W}'_k(\tau)$ are normalizations of the Wronskians
%$W_k(q)$ and $W'_k(q)$, respectively, so the previous congruence
%will hold if we can find constants $\lambda_k$, $\nu_k$,
%$g.c.d(\lambda_k,p)=g.c.d(\nu,k)=1$, such that
%$$ \lambda_k \frac{{W}'_k(\tau)}{{W}_k(\tau)} \equiv \nu_k \
%({\rm mod} p).$$
To prove (\ref{a}) we will use some row operations on the
determinant $W'_{L(k \Lambda_0)}$. Firstly, by the Leibnitz rule
$${\rm ch}_{k,i}^{(r)}(\tau)=\sum_{j=0}^{r} {r \choose j} \theta^{(j)}_{k,i}(\tau)
\left(\frac{1}{\eta(\tau)^3}\right)^{(r-j)}.$$
%We expand the numerator $\mathcal{W}'_k$ along the $(k+1)$-st row
%so
%$$\mathcal{W}'_k(z)=\sum_{i=1}^{k+1} {\rm ch}^{(k+1)}_{k,i}(z)M_i(z),$$
%where $M_i$ are the corresponding minors.
By applying a sequence of row operations, we first rewrite the
determinant $W'_{L(k \Lambda_0)}(\tau)$ so that in the $j$-th row
and $i$-th column we have
$$\frac{\theta^{(j)}_{k,i}(\tau)}{\eta(\tau)^3}+
a_{j}(\tau)\theta_{k,i}(\tau),$$ where $a_{j}(\tau)$ does not
depend on $i$. For instance, in the first row the entries are now
$$\theta_{k,i}(\tau)'
\frac{1}{\eta(\tau)^3}+\theta_{k,i}(\tau)\left(\frac{1}{\eta(\tau)^3}\right)',
\ i=1,...,k+1,$$ so
$$a_1(\tau)=\left(\frac{1}{\eta(\tau)^3}\right)'.$$ In the second
row the entries are
$$\frac{\theta^{''}_{k,i}(\tau)}{\eta(\tau)^3}+ \left(\left(\frac{1}{\eta(\tau)^3}\right)''-
2 \eta(\tau)^3 \left(\left(\frac{1}{\eta(\tau)^3}\right)'\right)^2
\right) \theta_{k,i}(\tau), \ i=1,...,k+1,$$ thus
$$a_2(\tau)=\left(\frac{1}{\eta(\tau)^3}\right)''- 2
\eta(\tau)^3 \left(\left(\frac{1}{\eta(\tau)^3}\right)'\right)^2
.$$ Clearly, $a_j(\tau)$ can be defined recursively. Let
$$A_j(q)=\eta(\tau)^3 a_j(\tau).$$
Consider an exponential generating function
$$\mathcal{A}(q)=\sum_{n=0}^\infty A_n(q) \frac{y^n}{n!}.$$
We claim that
\begin{equation}\label{geninv} \sum_{n=0}^\infty
A_n(q) \frac{y^n}{n!}=-\left(\frac{\eta(e^y q)}{\eta(q)}\right)^3.
\end{equation}
To see this notice that the exponential generating function for
{moments} of $\frac{1}{\eta(q)^3}$ is given by
\begin{equation} \label{moment}
\left(\frac{\eta(q)}{\eta(e^y q)}\right)^3,
\end{equation}
where the $j$-th moment of $\frac{1}{\eta^3(q)}$ is
$$\eta(q)^3 \left(\frac{1}{\eta(q)^3} \right)^{(j)}.$$
%\left(\eta(\tau)^{-3}\right)^{(j)}}{\eta(\tau)^{-3}}.$$
Now, the generating function $\mathcal{A}(q)$ is just $(-1)$ times
the reciprocal of the generating function for the moments
(\ref{moment}). This can be seen from the recursion formula for
$A_j(q)$'s. Now we have a closed expression for $a_j(\tau)$, so in
the $j$-th row the entries are
\begin{equation} \label{reduced}
\frac{\theta^{(j)}_{k,i}(\tau)}{\eta(\tau)^3}-j! \left([y^j]
\frac{\eta^3(e^y q)}{\eta^3(q)}\right)
 \frac{\theta_{k,i}(\tau)}{\eta^3(\tau)}, \ i=1,...,k+1.
\end{equation}
Let us focus at the $(k+1)$-st row. The formula (\ref{reduced})
together with Proposition \ref{jacmodp} and Lemma \ref{k1st} imply
\begin{eqnarray}
&& \frac{\theta^{(k+1)}_{k,i}(\tau)}{\eta^3(\tau)}+ a_{k+1}(\tau)
\theta_{k,i}(\tau) \nonumber \\
&& \equiv 4(k+2)^{-k-1}
\frac{\theta_{k,i}(\tau)}{\eta^3(\tau)}-2^{-3(k+1)}
\frac{\theta_{k,i}(\tau)}{\eta^3(\tau)} \ ({\rm mod} \ p).
\end{eqnarray}
Finally, the formula
$$(4(k+2))^{-k-1}-2^{-3(k+1)} \equiv 0 \ ({\rm mod} \ p)$$
implies
$$W'_k \equiv \frac{W'_k}{W_k} \equiv 0 \ ({\rm mod} \ p).$$
The proof now follows.
\end{proof}
\begin{conjecture} \label{ssp}
For every $k \in \mathbb{N}$, and $p=2k+3 \geq 5$ prime,
$$\mathcal{F}_{L(k \Lambda_0)}(q)  \equiv 1 \ ( {\rm mod} \ p).$$
\end{conjecture}
%\noindent {\bf I've tried a few obvious things (e.g., row
%operations, etc.), but it didn't work.}

\section{The zeros of $G(\mathcal{F}_{L(k \Lambda_0)},j)$}

Let us recall that every holomorphic modular form $f(\tau)$ can be
uniquely expressed as \be \label{rep} f(\tau)=\Delta^t(\tau)
E_4^\delta(\tau) E_6^\epsilon(\tau) {G}(f,j(\tau)), \ee where
$$j(\tau)=\frac{1728 E_4^3(\tau)}{E_4^3(\tau)-E_6^2(\tau)},$$
${G}(f,j)$ is a polynomial of degree $\leq t$ and
$$k=12t+4 \delta+6 \epsilon,$$
where $0 \leq \delta \leq 2$ and $0 \leq \epsilon \leq 1$. It is
known that the $j$-function defines a one-to-one map from the arc
$[e^{2 \pi i/3},e^{\pi i/2}]$ onto the interval $[0,1728]$. As in
\cite{MMO}, based on extensive computations we conjecture
\begin{conjecture} \label{two-conj}
For every $k \neq 2i^2-2$, $i \geq 2$, the zeros of
$G(\mathcal{F}_{L(k \Lambda_0)},j)$ are simple and inside the
interval $[0,1728]$.
\end{conjecture}
\noindent Here is a sample of $j$-zeros for $1 \leq k \leq 11$
which clearly supports our conjecture.
\begin{center} \label{table1}
\begin{tabular}{|c|c|c|c|c|c|}
\hline
$k$ & $\epsilon$ & $\delta$    & $G(\mathcal{F}_{L(k \Lambda_0)},j)$ & zeros of ${G}(\mathcal{F}_{L(k \Lambda_0)},j)$ \\
\hline \hline
$1$ & $1$   & $0$ &   $1$ & $$  \\
\hline
$2$ & $0$  &  $1$ &  $1$ &  $$   \\
\hline
$3$  & $2$   &  $0$ & $1$  &    $$  \\
\hline
$4$ & $1$   &  $1$ &  $1$  &   $$ \\
\hline
$5$ & $0$   &  $0$ &  $j-\frac{1302528}{1075}$ & $1211.653954$  \\
\hline
$6$ & $2$   &  $1$ &  $0$ & $-$  \\
\hline
$7$ & $1$   &  $0$ &  $j-\frac{787021824}{587489}$ & $1339.636698$  \\
\hline
$8$ & $0$   &  $1$ &   $j-\frac{8696400}{20119}$ & $432.2481237$  \\
\hline
$9$ & $2$   &  $0$ &  $j-\frac{1381580800}{10776887}$ & $1281.987070$  \\
\hline
$10$ & $1$   & $1$  &  $j-\frac{956352}{2021}$ & $473.2073231$  \\
\hline $11$ & $0$ & $0$  &
$j^2-\frac{20462710947840}{13928908741}j+\frac{1908473415598080}{13928908741}$
&
$100.0843760, 1368.997756$  \\
\hline
\end{tabular}
\end{center}

\section{Final remarks}

Notice that (\ref{hasse}) would follow from the congruence
$$W'_{L(k \Lambda_0)}(q) \equiv  h W_{L(k \Lambda_0)}(q) \ ({\rm
mod} \ p^2),$$ for some $h$.

The method in \cite{R} is peculiar to Eisenstein series and it
does not apply directly to $F_{L(k \Lambda_0)}$. So in order to
probe Conjecture \ref{myfirst}, E. Mortenson computed the
$j$-zeros of $\mathcal{F}_{L(k \Lambda_0)}$ for all $1 \leq k \leq
22$. From his data we observed certain interlacing properties of
zeros, which indicates a possibility that $G(\mathcal{F}_{L(k
\Lambda_0)},j)$, $k \in \mathbb{N}$ forms an orthogonal polynomial
sequence (cf. \cite{KZ1}). An unpleasant feature of $F_{L(k
\Lambda_0)}$ is rather irregular pattern of vanishing which occurs
for $k=2i^2-2$. This, and some other clues, makes it hard to
believe that a simple recursive formula for $F_{L(k \Lambda_0)}$
will settle down the conjecture.
We should mention here that recursion formulas do arise naturally
in the context of $L(\Lambda_0)^{\otimes^k}$-modules (see
\cite{M3}). But as we know $L(k \Lambda_0)$-modules and
$L(\Lambda_0)^{\otimes^k}$-modules are related in a nontrivial way
via the unitary minimal models.

\end{document}